\documentstyle[a4,12pt]{article}
\begin{document}
\newcounter{tn}
\setcounter{tn}{0}
\newtheorem{theo}{Theorem}[section]
\newtheorem{dfn}[theo]{Definition}
\newtheorem{prop}[theo]{Proposition}
\newtheorem{lem}[theo]{Lemma}
\newtheorem{coro}[theo]{Corollary}
\newenvironment{mth}{\medskip \par \noindent
 {\bf Main Theorem} \it}{\rm \smallskip \par \noindent}
\newcommand{\bp}{\noindent{\em Proof. }}
\newcommand{\ep}{\hfill$\Box$\par\medskip}
\newcommand{\nop}{\hfill$\Box$\par}
\def\aut{\mbox{\rm Aut}}
\newcounter{szam}
\newcounter{mul}
\hyphenation{a-u-to-mor-phism  a-u-to-mor-phisms} 
%
%
%
%
%
%
%
%
%
\title{Collineation groups of the smallest Bol
$3$-nets\footnote{This paper is accepted for publication in {\em
Note di Matematica}, vol. 14 (1994), n.1.}}
\author{G\'abor Nagy\thanks{The author wishes to thank the
opportunity of the participation in the {\it Bolyai Network}
TEMPUS JEP project.}} 
\maketitle
\begin{abstract}
Some associativity properties of a loop can be interpreted as
certain closure configuration of the corresponding $3$-net. It
was known that 
the smallest non-associative loops with the so called left Bol
property have order 8. In this paper, we determine the direction
preserving collineation groups of the $3$-nets belonging to these
smallest Bol loops. For that, we prove some new results
concerning collineations of $3$-nets and autotopisms of loops.
\end{abstract}
%
%
%
%
\section{Introduction}
The theory of $3$-nets and loops is a discipline
which takes ist roots from geometry, algebra and combinatorics.
In geometry, it aroses from the analysis of web structures; in
algebra, from non-associative products; and in combinatorics,
from Latin Squares. The subject has grown by relating these
aspects to an increasing variety of new fields, and is now
developing into a thriving branch of mathematics.

The reader of this paper is referred to {\sc H.O. Pflugfelder}'s
book \cite{Pflug} for the definition
of the {\em loop}, the {\em nuclei} and the {\em left and right
translations} of a given loop (\cite{Pflug}, p. 3., 5., 16.).
The concept of a $3$-net can be formulated in many equivalent
ways, for us, if a $3$-net is coordinatized by the loop
$(L,\cdot)$, then the transversal lines are sets of points of the
form $\{ (x, y) \in L \times L : x \cdot y = \mbox{constant}\}$.
We define the {\em section} $S(L)$ of the loop $(L,\cdot)$ by
$S(L) = \{\lambda_x : x\in L\}$, and the {\em group generated by
the left translations} of $(L,\cdot)$ is $G(L) = \langle S(L)
\rangle$. 

There was not much work made in the direction of collineation
groups of $3$-nets, few exceptions are \cite{Barl}, \cite{Boni},
\cite{Saeli}. On the other hand, the class of Bol loops are also
connected with configurations of the projective plate, cf.
\cite{Kor}, \cite{Kall}. 

We have the
\begin{dfn}
A bijection $\gamma$ on $L$ is called a right (left) {\em
pseudo-auto\-mor\-phism} of a loop $(L,\cdot)$ if there exists at
least one element $c \in L$ such that 
\[x^\gamma \cdot (y^\gamma \cdot c) = (xy)^\gamma \cdot c \]
$((c \cdot x^\gamma) \cdot y^\gamma = c \cdot
(xy)^\gamma)$ for all $x,y \in L$. The element $c$ is then
called a {\em companion} of $\gamma$.
\end{dfn}
This concept will play a significant role in the examination
of the col\-line\-ations of a $3$-net. The next proposition gives a
list of the most important properties of pseudo-automorphisms.
\begin{prop} \label{pseu} Let $\gamma$ be a right
pseudo-automorphism of a loop $(L,\cdot)$ with a companion $c$.
Then it has the following properties:
\setcounter{szam}{0}
\begin{list}{(\roman{szam})}{\usecounter{szam}}
\item $1^\gamma = 1$.
$\gamma^{-1}$ is a right pseudo-automorphism with a companion
$d$ such that $c^{\gamma^{-1}} \cdot d = 1$.
If $\gamma^\prime$ is another pseudo-automorphism of $L$
with a companion $c^\prime$ then $\gamma \gamma^\prime$ is a
pseudo-automorphism of $L$ with the companion $c^\gamma \cdot
c^\prime$. 
The right pseudo-automorphisms of a given loop form a group.
\item It holds $n^\gamma\cdot x^\gamma = (nx)^\gamma$ for any
$n\in N_\lambda$ and $x \in L$. The permutation $\gamma$ induces
an automorphism of the left nucleus of $(L,\cdot)$.
\item It holds $x^\gamma\cdot n^\gamma = (xn)^\gamma$ for any
$n\in N_\mu$ and $x\in L$. The permutation $\gamma$ induces an
automorphism of the middle nucleus of $(L,\cdot)$.
\end{list} 
Analogous statements are valid for left pseudo-automorphisms.
\end{prop}
\bp The statements of $(i)$ are easy to
prove with a few calculation (see \cite{Pflug}, p. 74-76.). In
\cite{Pflug}, p. 77. we find the points $(ii)$ and $(iii)$, but
they are stated in an incorrect way, and also the proof contains
some technical faults, so we give now a complete proof.

Let $\gamma$ be a right pseudo-automorphism with a companion $c$
of $(L, \cdot)$ with unit element $1$. Then
\[(n x) y = n(x y) \> \> \> \mbox{ for all } \> \> \> x, y \in L
\> \> \> \mbox{ and } \> \> \> n\in N_\lambda.\]
With $\langle \gamma, \gamma \rho_c, \gamma \rho_c \rangle$
applied to this, we have 
\[(n x)^\gamma \cdot (y^\gamma \cdot c) = n^\gamma \cdot
(x^\gamma (y^\gamma \cdot c)).\]
Replace $y^\gamma \cdot c$ by $z$ to obtain
\[(n x)^\gamma \cdot z = n^\gamma (x^\gamma \cdot z),\]
or with $z = 1$,
\[(n x)^\gamma = n^\gamma \cdot x^\gamma.\]
Comparing the two last equations, one has
\[(n^\gamma \cdot x^\gamma) \cdot z = n^\gamma \cdot (x^\gamma
\cdot z),\]
which implies that $n^\gamma \in N_\lambda$, thus $(ii)$ is
shown. Similarly, we can prove for $(iii)$ as well. \ep

A {\em collineation} of a 3-net is a one-to-one map on the point
set, such that two points are collinear if and only if their
images are collinear. We speak about a {\em direction preserving
collineation} if for any line of the $3$-net, the line and its
image belong to the same paralell class. If we want to study the
collineation group of a 
3-net, it could be useful to be able to write down a
collineation map with algebraic tools, using the
coordination of the 3-net.
\begin{prop} \label{dircoll}
Let $\cal N$ be a $3$-net coordinatized by the loop $(L,\cdot)$ and
$\gamma$ be a permutation of the point set of $\cal N$ such that
$(x,y)^\gamma = (x^{\alpha_y},y^{\beta_x})$, with $\alpha_y,
\beta_x$ permutations of $L$. Then $\gamma$ is a direction
preserving collineation if and only if
\setcounter{szam}{0}
\begin{list}{(\roman{szam})}{\usecounter{szam}}
\item $\alpha_y$ and $\beta_x$ do not depend on $y$ and $x$,
respectively, this means that we can write $\alpha=\alpha_y$ and
$\beta=\beta_x$. 
\item For every $x,y \in L$ holds
\[x^\alpha \cdot y^\beta = 1^\alpha \cdot (xy)^\beta.\]
\end{list} \end{prop}
\bp Suppose that $\gamma$ is a direction preserving
collineation. Then it preserves the vertical lines, so the first
coordinate of the image of a point depends only on the first
coordinate of the point. The same idea holds for the horizontal
lines with the second coordinates, so $(i)$ follows. Moreover,
because of the definition of the 3-net, the points $(x,y)$ and 
$(1,xy)$ lie on the same transversal line. Then so do
$(x^\alpha,y^\beta)$ and $(1^\alpha,(xy)^\beta)$, i. e. $(ii)$
holds.

Conversely, let $\gamma$ be such that $(x,y)^\gamma =
(x^\alpha,y^\beta)$ and $\alpha$ and $\beta$ satisfies $(ii)$.
It is easy to see that $\gamma$ maps vertical and horizontal
lines to vertical and horizontal lines, respectively. And it
preserves the third direction as well, because the points
$(x,y)$ and $(x^\prime, y^\prime)$ are on a transversal line iff
$x \cdot y = x^\prime \cdot y^\prime$. Then 
\[x^\alpha \cdot y^\beta = 1^\alpha \cdot (xy)^\beta = 
1^\alpha \cdot (x^\prime y^\prime)^\beta = 
{x^\prime}^\alpha \cdot {y^\prime}^\beta,\]
thus the images of the two points lay on a transversal line,
too. \ep

As an application of this Proposition, we can determine certain
stabiliser subgroups of the direction preserving collineation
groups.

Firstly, let us suppose that $(\alpha,\beta)$ is a collineation
fixing the horizontal line through the origin, that means that
$1^\alpha = c$ and $1^\beta=1$. Applying $(\alpha,\beta)$ for
the point $(x,1)$, from (ii), we get $x^\alpha = x^\alpha \cdot
1^\beta = 1^\alpha \cdot x^\beta = c\cdot x^\beta$. Now, for
every $x,y \in L$ we have 
\[(c\cdot x^\beta)\cdot y^\beta = x^\alpha \cdot y^\beta =
c\cdot (xy)^\beta,\]
thus, a collineation which fixes the horizontal line through the
origin determines a left pseudo-automorphism with the companion
$c = 1^\alpha$. Similarly, we can prove the link between
stabilisers of the vertical line and right pseudo-automorphism.
Finally, if $(\alpha, \beta)$ fixes the origin, $1^\alpha =
1^\beta = 1$, than $\alpha=\beta$ are automorphisms of the
coordinate loop.

Moreover, let $P$ denote the orbit of the origin in the action
of the group of the direction preserving collineations, and
$l_h$ and $l_v$ the horizontal and vertical lines through the
origin, respectively. Then we have the natural bijections $l_h
\cap P \rightarrow C_\lambda$ and $l_v \cap P \rightarrow
C_\rho$, where $C_\lambda$ and $C_\rho$ are the sets of those
elements of $L$ which can occur as left and right companions,
respectively. 

It is known that in a given $3$-net, the coordinate loop
depends on the choice of the origin point. The following
proposition gives a sufficient and necessary conditions for
that, that the coordinate loops with the origins $P$ and $Q$ be
isomorphic. 
\begin{prop} \label{isom}
Let $P$ and $Q$ be points of the $3$-net $\cal N$ and let
$(L,\cdot)$ and $(L,\circ)$ be the coordinate loops with the origin
$P$ and $Q$, respectively. These two loops are isomorphic if and
only if there is a direction preserving collineation $\gamma$ of
the net, mapping $P$ onto $Q=P^\gamma$.
\end{prop}
\bp Clearly, if the collineation $\gamma$ exists, then it
induces a ono-to-one map between the point sets of the
horizontal lines through 
$P$ and $Q$, which leads to a permutation of $L$, which is going
to be an isomorphism. Conversely, let $U$ be an isomorphism and
$R$ be a point of $\cal N$ with coordinates $(x,y)$ in the
coordinate system belonging to $P$. Define the image of $R$ by
the point $R^\prime$, which has the coordinates $(x^U, y^U)$ in
the coordinate system belonging to $Q$. Obviously, the map
$\gamma : R \rightarrow R^\prime$ is a direction preserving
collineation of $\cal N$.
\ep
%
%
%
%
\section{Special loop classes}
If we have an algebraic structure with a non-associative
multiplication, then one can ask questions about the necessary
and sufficient conditions for special algebraic identities which
interpret weak associative properties. In the most cases, such
identities lead to geometrical configurations in the
coordinatized $3$-net. 

We say that a loop satisfies the {\em left inverse property}
(L.I.P.), if the inverse of any left translation is contained in
the section. With other words, this means that for any element
$x$ of the loop there is a unique element $x^{-1}$ such that
$\lambda_x^{-1} = \lambda_{x^{-1}}$. Let us define the map $J: L
\rightarrow L$ by $J: x \mapsto x^{-1}$.
\begin{prop} Let $(L,\cdot)$ be a loop satisfying the left
inverse property, and let $N_\lambda$ and $N_\mu$ be its left
and middle nuclei, respectively. Then we have $N_\lambda =
N_\mu$. 
\end{prop}
\bp Suppose $n\in N_\lambda$, then also $n^{-1} \in N_\lambda$.
This is equivalent to the fact that for all $x\in L$, we have
$\lambda_x \lambda_{n^{-1}} \in S(L)$, or, we can write
$\lambda_{x^{-1}} \lambda_{n^{-1}} \in S(L)$. Using the left
inverse property, we obtain that for all $x\in L$
\[\lambda_n \lambda_x = \lambda_{n^{-1}}^{-1}
\lambda_{x^{-1}}^{-1} = (\lambda_{x^{-1}} \lambda_{n^{-1}})^{-1}
\in S(L).\]
And this holds if and only if $n\in N_\mu$ as well. \ep

An important class of loops having the L.I.P. is the so called
{\em Bol loops}, which can be characterized by an algebraic
identity, or, equivalently, by a closure configuration in the
$3$-net. {\sc G. Bol}'s configuration and the corresponding
algebraic identity have always played an important role in the
development of the geometric theory of $3$-nets and in the
algebraic theory of loops. 

\begin{dfn} A loop $(L,\cdot)$ is said to be a {\em (left) Bol
loop} if the identity 
\[x\cdot (y \cdot (x \cdot z)) = (x\cdot (y\cdot x)) \cdot z\] 
holds true for all $x,y,z \in L$. 
\end{dfn}
Or, we can formulate this, saying that $(L,\cdot)$ is a
(left) Bol loop if for all $x,y \in L$ we have
\[\lambda_x \lambda_y \lambda_x \in S(L).\]

In a Bol $3$-net, the Bol configuration defines a collineation
which is a reflection with a vertical axis.
In the paper \cite{Funk}, the subgroup of the direction
preserving collinaetion group $\Gamma(\cal N)$
containing the products of these Bol reflections with an
even number of factors was considered. Let us denote this
subgroup by $N(\cal N)$. 
\begin{prop}
The set of Bol reflections of a Bol $3$-net is a normal set in
the full group of collineations.
\end{prop}
\bp Let us observe that the
paralell class of the vertical lines is special, because the Bol
reflection by them is equivalent with the left Bol condition,
hence the full group of collineations of the $3$-net leaves the
set of vertical lines fixed. Let $\sigma_m$ be the Bol
reflection with the vertical line $(x = m)$ as axis. Then, the
unique collineation which fixes the axis pointwise and
interchanges the horizontal and transversal directions is
$\sigma_m$. Hence, if we conjugate $\sigma_m$ with an arbitrary
collineation $\gamma$, we get the Bol reflection
$\sigma_{m^\prime}$, where the vertical line $(x = m^\prime)$ is
the image of the axis $(x = m)$ by $\gamma$. \ep
It follows that $N(\cal N)$ is normal in the group of direction
preserving col\-line\-ations. They have shown that there is an
endomorphism $\Phi: N({\cal N}) \rightarrow G(L)$, where $L$ is
a coordinate loop of $\cal N$ 
and $G(L)$ is the group generated by the left translations of
$L$. Moreover, the kernel of $\Phi$ is isomorphic to a subgroup
of the left nucleus of $L$. We also know that $N(\cal N)$ acts
transitively on the set of the horizontal lines. In the paper, a
new coordinatization of the 3-net was introduced. The origin was
the same, and if a point had coordinates $(x, y)$ in the old
system, then in the new one it is $(\xi, \eta)$, with $\xi = y$ and
$\eta = xy$. It was shown, that this coordinate system gives the
same multiplication on the point set of the horizontal line
through the origin. Then, $N(\cal N)$ is generated by the
collineations of the form 
\[(\xi, \eta) \mapsto (\bar\xi, \bar\eta) = (a^{-1}\xi,
a\eta),\] 
with $a\in L$. Let us translate this in the old coordinate
system. Define $(\bar x, \bar y)$ as the image of the point $(x,
y)$ by the above collineation. Then we have
\[\bar y = \bar\xi = a^{-1}\xi = a^{-1} y = y^{\lambda_a^{-1}},\]
and
\[\bar x = \bar\eta / \bar\xi = (a \eta) / (a^{-1}\xi) = (a
\cdot xy) / (a^{-1} y) = a \cdot x a = x^{\rho_a \lambda_a}.\] 
The last step follows from the fact, that this is a direction
preserving col\-line\-ation, thus the $x$-coordinate of the image
point does not depend on the $y$-coordinate, and so we can
choose $y=a$. We got that the generators of $N(\cal N)$ are of
the form $(\rho_a \lambda_a, \lambda_a^{-1})$, where $a\in L$. 
The map $\Phi$ is a projection of the second permutation.

A useful characterization of the Bol loops are provided by the
following 
\begin{prop} \label{karclo}
A loop $(L,\cdot)$ is a Bol loop if and only if for
any $x\in L$, the pair of permutations $(\rho_x \lambda_x,
\lambda_x^{-1})$ is a direction preserving collineation of the
$3$-net, coordinatized by the loop $(L,\cdot)$.
\end{prop}
\bp By the above, it is obvious that a Bol $3$-net has
collineations of the form $(\rho_x \lambda_x, \lambda_x^{-1})$.
Conversely, if we assume that $(\rho_x \lambda_x,
\lambda_x^{-1})$ is a collineation, than using \ref{dircoll}
for the point $(u, xv)$, we have
\[(x \cdot ux)\cdot v = u^{\rho_x \lambda_x} \cdot
v^{\lambda_x^{-1}} = 1^{\rho_x \lambda_x} \cdot
(u \cdot xv)^{\lambda_x^{-1}} = x ( u\cdot xv),\]
thus, the Bol identity holds for every $u,v \in L$.
\ep

Concerning the kernel of $\Phi$ we have the
\begin{prop} \label{kernel}
Let $\cal N$ be a Bol 3-net coordinatized by the Bol loop
$(L,\cdot)$. Suppose that $J$ is a right pseudo-auto\-mor\-phism
of $L$. Then the kernel of $\Phi$ is trivial.
\end{prop}
\bp The generators of $N(\cal N)$ are in the form
$(\rho_x \lambda_x, \lambda_x^{-1})$. Then one can easily prove
by induction on $k$ that the
elements of $N(\cal N)$ are in the form 
\[(\rho_{m_1}\lambda_{m_1} \cdots \rho_{m_k}\lambda_{m_k},
\lambda_{m_1}^{-1} \cdots \lambda_{m_k}^{-1}) =
(\rho_u\lambda_{m_1} \cdots \lambda_{m_k}, 
\lambda_{m_1}^{-1} \cdots \lambda_{m_k}^{-1}), \]
with $u = m_1(\cdots(m_{k-2}\cdot m_{k-1} m_k)\cdots ) =
1^{\lambda_{m_k} \cdots \lambda_{m_1}}$. 
Obviously, 
\[\Phi: (\rho_{m_1}\lambda_{m_1} \cdots \rho_{m_k}\lambda_{m_k},
\lambda_{m_1}^{-1} \cdots \lambda_{m_k}^{-1})
\mapsto \lambda_{m_1}^{-1} \cdots \lambda_{m_k}^{-1}.\]
Suppose that
$(\alpha,\beta) \in Ker \Phi$, then $\beta = id$.
By \ref{dircoll}, $\alpha$ must be of the form $\lambda_a$ with
$a=1^\alpha$. Since $id = \beta =\lambda_{m_1}^{-1} \cdots
\lambda_{m_k}^{-1}$, we have $m_k^{-1}(\cdots(m_3^{-1}\cdot
m_2^{-1} m_1^{-1})\cdots ) =1$. 
Let $c$ be a companion element of $J$. Then $u=1^{\beta^{-1}} =
1$ and
\[a c = c^\alpha = m_k(\cdots(m_3(m_2 \cdot m_1 c))\cdots )=\] 
\[= m_k(\cdots(m_3((m_2^{-1} m_1^{-1})^{-1}\cdot c))\cdots )=
m_k(\cdots((m_3^{-1}\cdot m_2^{-1} m_1^{-1})^{-1} c)\cdots ) =\]
\[= \cdots = (m_k^{-1}(\cdots(m_3^{-1}\cdot m_2^{-1}
m_1^{-1}))\cdots )^{-1} \cdot c = 1^{-1} \cdot c = c,\]
thus $a=1$ and $\alpha=id$. \ep
{\it Remark.} Actually, we have proven now a little more than
the theorem. We have just shown that if $(\alpha,\beta) \in
N({\cal N})$ and $1^\beta = 1$, then $\alpha$ fixes those elements
of the loop which are companions of the right
pseudo-auto\-mor\-phism $J$.

We are interested in another class of loops as well, namely the
class of loops $(L,\cdot)$ for 
which $S(L)$ is closed under conjugation in the symmetric group
on $L$ (cf. \cite{NSt}). More precisely, we make the following
\begin{dfn}
We say that the loop $(L,\cdot)$ is a {\em left conjugacy closed
loop} if 
\[\lambda_x^{-1} \lambda_y \lambda_x \in S(L)\]
holds for all $x,y \in L$.
\end{dfn}
This class of loops can be characterized by right
pseudo-automorphisms of a given form.
\begin{prop}
A loop $(L,\cdot)$ is left conjugacy closed if and only if for
all $x\in L$, the permutation $\lambda_x \rho_x^{-1}$ is a right
pseudo-automorphism with the companion element $x$.
\end{prop}
\bp We have $u^{\lambda_x \rho_x^{-1}} \cdot (v^{\lambda_x
\rho_x^{-1}} \cdot x) = (x u) / x \cdot (x v)$ and $(u
v)^{\lambda_x \rho_x^{-1}} \cdot x = x \cdot uv = x\cdot ( u
\cdot ( x^{-1} \cdot x v ) )$. Setting $\bar{v} = xv$, we got
that the permutation $\lambda_x \rho_x^{-1}$ is a right
pseudo-automorphism with companion $x$ exactly if $(x u) / x
\cdot \bar{v} = x\cdot ( u \cdot ( x^{-1} \bar{v}) )$ for all
$u,\bar{v},x \in L$, i.e. if $\lambda_x^{-1} \lambda_u \lambda_x
\in S(L)$ for all $u,x \in L$. \ep 
\begin{coro}
Let $\cal N$ be a 3-net coordinatized by a left conjugacy closed
loop and consider the action of the collineation group on the
point set of $\cal N$. Then the orbit of the origin is a union
of vertical lines. \end{coro}
\bp By \ref{karclo}, the orbit contains the vertical line through
the origin. Now, suppose that it contains a point $P$ outside of
this line. Then we coordinatize with this point, and by
\ref{isom}, the coordinate loop is going to be isomorphic to the
original one, hence it is left conjugacy closed, and so the orbit
contains again the whole vertical line through $P$. \ep
%
%
%
%
\section{Automorphisms of the smallest Bol loops}

The smallest Bol loop which is not a group must have
eight elements. These loops form two isotopy classes;
representatives of the isotopy classes are given by the
following (cf. \cite{Burn}).
\setcounter{szam}{0}
\begin{list}{\arabic{szam}.)}{\usecounter{szam}}
\item The loop $B_1$ has as elements $e$, $f$, $f^2$,
$f^3$, $g$, $fg$, $f^2g$, $f^3g$ such that the multiplication is
defined by the following rules
\setcounter{mul}{0}
\begin{list}{(\roman{mul})}{\usecounter{mul}}
\item $f^i(f^j g) = f^{i+j} g$, 
\item $(f^i g)f^j = f^{i+j} g$,
\item $(f^i g)(f^j g) = f^{2-i+j}$, $i,j \in Z_4$,
\end{list} 
with $i,j \in Z_4$.
\item The loop $B_2$ has as elements $e$, $f$, $f^2$,
$f^3$, $g$, $fg$, $f^2g$, $f^3g$ such that the multiplication is
defined by the following rules
\setcounter{mul}{0}
\begin{list}{(\roman{mul})}{\usecounter{mul}}
\item $f^i(f^j g) = f^{i+j}g$, 
\item $(f^i g)f^j = f^{i+j+2ij} g$,
\item $(f^i g)(f^j g) = f^{i+j+2ij}$,
\end{list} 
with $i,j \in Z_4$.
\end{list} 
The unit element is $e$, and for both cases, the left nucleus
contains the elements $e$ and $f^2$, it has order 2. The right
nuclei have order 4, but in the first case, $N_\rho = \{ e, f,
f^2, f^3\} \cong C_4$, and in the second case $N_\rho = \{ e, g,
f^2, g f^2 \} \cong C_2 \times C_2$.

In this chapter, we are going to examine the automorphism group
of the loop $B_1$. \ref{rovb2} will state that almost the
same statements hold for $B_2$, as well. For this theorem,
we will not give the whole proof, just mention the most
important differences, compared to the proofs of the first part.

\begin{prop} \label{b11}
The maps $\phi_1: f \mapsto f, g \mapsto fg$ and $\phi_2: f
\mapsto f^{-1}, g \mapsto g$ can be extended to an automorphism
of $B_1$. The permutation $\phi_1$ has order $4$, the permutation
$\phi_2$ has order $2$. Let $\phi_3 = \phi_1 \circ \phi_2$. Then
$\phi_3$ has order $2$. 
\end{prop}
\bp Let us define $\hat g$ by $\hat g = fg$. We will show that
the same multiplication rules hold as above if we write $\hat g$
in the place of $g$.
\setcounter{mul}{0}
\begin{list}{(\roman{mul})}{\usecounter{mul}}
\item $f^i(f^j \hat g) = f^i(f^j (fg)) = f^i(f^{j+1} g) = f^{i+j+1} g =
f^{i+j} (fg) = f^{i+j} \hat g;$
\item $(f^i \hat g)f^j = (f^{i+1} g)f^j = f^{i+j+1} g = f^{i+j} \hat g;$
\item $(f^i \hat g)(f^j \hat g) = (f^{i+1} g)(f^{j+1} g) = f^{2-i+j}.$
\end{list} 
By this we have proven that $\phi_1$ can be extended to an
automorphism of $B_1$. From now on, we will call the
automorphism itself $\phi_1$. 
In a very similar way we can prove for $\phi_2$, and we will use
this notation also for the automorphism. 

The order of $\phi_1$ follows from the first multiplication
rule, namely $g^{\phi_1^2} = f^\phi_1 \cdot g^\phi_1 = f(fg) =
f^2 g$, and by induction, $g^{\phi_1^i} = f^i g$. Clearly,
$\phi_2$ is an involution. The action of $\phi_3$ on the
generators is $f \mapsto \hat f = f^{-1}$, $g \mapsto \hat g = fg$.
Applying $\phi_3$ once again, we get $\hat f \mapsto {\hat f}^{-1} =
(f^{-1})^{-1} = f$, $\hat g \mapsto \hat f \hat g = f^{-1} (fg) = g$, so
$\phi_3^2 = id$. \ep
\begin{theo}
In the loop $B_1$, each left pseudo-automorphism is an
automorphism, only the elements of the left nucleus occur as 
left companions. And the group of automorphisms of $B_1$ is
isomorphic to the dihedral group $D_8$ of order $8$.
\end{theo}
\bp By \ref{b11}, the subgroup $\langle
\phi_1,\phi_2 \rangle \cong D_8$ generated by the permutations
$\phi_1$ and $\phi_2$ has the relations $\phi_1^4 = id$,
$\phi_2^2 = id$ and $(\phi_1 \circ \phi_2)^2 = id$ on its
generators, thus this subgroup is isomorphic to the dihedral
group of order 8. 

On the other hand, we have an upper bound for the order of the
group of left pseudo-automorphisms. By \ref{pseu} $(ii)$, a
left pseudo-automorphism induces an automorphism on $N_\mu =
N_\lambda$, so it must fix $e$ and $f^2$. Hence we have two
choices for the image of $f$, itself or $f^3$. For the image of
$g$ we get at most 4 possibilities, since $g$ is not contained
in $N_\rho$, so neither its image. Now, we know by \ref{pseu}
$(ii)$ as well, that the images
of $f$ and $g$ determine uniquely the pseudo-automorphism. 
For their images, we have altogether at most
$2\cdot4=8$ choices, hence $B_1$ has at least 8 automorphisms. These
are also pseudo-automorphisms (one can have the unit element as
companion), so there is no more pseudo-automorphism. This also
means, that there are exactly 8 automorphisms, hence their group
is the dihedral group $D_8$. \ep

Now, let us consider the case of the loop $B_2$.
\begin{theo} \label{rovb2} \setcounter{szam}{0}
\hfill
\begin{list}{(\roman{szam})}{\usecounter{szam}}
\item The maps $\psi_1: f \mapsto fg, g \mapsto g$ and $\psi_2: f
\mapsto f, g \mapsto f^2g$ can be extended to an automorphism
of $B_2$. The permutations $\psi_1$ and $\psi_2$ have order $2$.
Let be $\psi_3 = \psi_1 \circ \psi_2$. 
Then $\psi_3 \neq \psi_2 \circ \psi_1$ has order $4$.
\item In the loop $B_2$, each left pseudo-automorphisms are
automorphisms, only the elements of the left nucleus occur as
left companions, and the group of automorphisms of $B_2$ is
isomorphic to $D_8$.
\end{list} \end{theo}
\bp As in the first part of the chapter, we could show (i) using
the multiplication rules of $B_2$. In the prove of $(ii)$, we
use that the 
left nucleus has 2 elements, the right nucleus has 4 elements.
In this case, it is not true that $N_\rho$ contains the powers
of $f$, but $N_\rho = \{e, g, f^2, f^2 g\}$ and $N_\lambda = \{e,
f^2\}$. Hence, we can show $(ii)$, as well. \ep

In the rest of this chapter, we prove two propositions which
give some extra information on the non-associative
multiplication structure of these smallest Bol loops.
\begin{prop} \label{stab}
Let be $a,b,c \in L$ such that $1^{\lambda_a \lambda_b
\lambda_c} = 1$, where $L = B_1$ or $L=B_2$. Then, $\lambda_a
\lambda_b \lambda_c = id$ iff at least one of the elements
$a,b,c$ belong to the left nucleus $N_\lambda$.
\end{prop}
\bp It is obvious that for any Bol loop, if the product of 2
left translations fixes the unit element $1$ of the loop, then
the product is the identity permutation. Consider the product of
permutations $\lambda_a \lambda_b \lambda_c$. 
If one of $a$, $b$ or $c$ belongs to $N_\lambda = N_\mu$ then
the $\lambda_a \lambda_b \lambda_c$ can be written as the
product of 2 left translations, and so $\lambda_a 
\lambda_b \lambda_c = id$. 

By $c \cdot ba = 1$, at least one of these elements belongs to
$N_\rho$, since $N_\rho$ has index 2 in $L$.
Assume that at least 2 of them belongs to $N_\rho$. Then
all of them are in $N_\rho$, and at least one is in $N_\lambda
\leq N_\rho$, since $[N_\rho : N_\lambda] = 2$. Thus, in this
case, $\lambda_a \lambda_b \lambda_c = id$. This also means that
for any $n,m \in N_\rho$, $\lambda_n \lambda_m \lambda_{mn}^{-1} =
id$, i.e. $\lambda_n \lambda_m = \lambda_{mn}$.

Notice that there exist elements $a^\prime, b^\prime, c^\prime
\in L$, such that $c^\prime \cdot b^\prime a^\prime = 1$ and
$\lambda_{a^\prime} \lambda_{b^\prime} \lambda_{c^\prime} \neq
id$, since $L$ is not a group. By the paragraph above, exactly
one element of $a^\prime, b^\prime, c^\prime$ must be in $N_\rho$.
Now, suppose that exactly one of the elements $a,b,c$ is in
$N_\rho$, but none of them is in $N_\lambda$. The permutation
$\lambda_a \lambda_b \lambda_c$ is always a left
pseudo-auto\-mor\-phism (cf. \cite{Funk}), thus 
for $L$ it is an automorphism, we denote it by $\tau$. For any
automorphism $\alpha$ and $x\in L$ we have $\alpha^{-1}
\lambda_x \alpha = \lambda_{x^\alpha}$, and so
\[\lambda_a \lambda_b \lambda_c = \tau = \lambda_a^{-1} \tau
\lambda_{a^\tau} = \lambda_b \lambda_c \lambda_{a^\tau}.\]
This means that we may assume w.l.o.g. that $a, a^\prime \in
N_\rho$. We claim that there exists exactly one automorphism
$\alpha$ of $L$ such that $a^\alpha = a^\prime$ $b^\alpha =
b^\prime$ and $c^\alpha = c^\prime$. We know
that $N_\rho \setminus N_\lambda$ and $L \setminus N_\rho$ are
orbits in the action of the automorphisms on the carrier set
$L$. Moreover, we have exactly $8 = |N_\rho \setminus N_\lambda|
\cdot |L \setminus N_\rho|$ automorphisms, so there is exactly
one $\alpha$ with $a^\alpha = a^\prime$ and $b^\alpha =
b^\prime$, and than $c^\alpha = c^\prime$ since $c \cdot ba =
c^\prime \cdot b^\prime a^\prime = 1$ and $\alpha$ is
automorphism. We got
\[id \neq \lambda_{a^\prime} \lambda_{b^\prime} \lambda_{c^\prime} =
\lambda_{a^\alpha} \lambda_{b^\alpha} \lambda_{c^\alpha} =
\alpha^{-1} \lambda_{a} \lambda_{b} \lambda_{c} \alpha,\]
thus $\lambda_{a} \lambda_{b} \lambda_{c} \neq id$. 
\ep
{\it Remark.} If $c \cdot ba =1$ then the automorphism $\lambda
_a \lambda_b \lambda_c$ fixes the elements of the right nucleus
or with other words, it commutes with the left translations by a
right nucleus element. 

\begin{prop} \label{nuc}
Let the loop $L$ be equal to the loop $B_1$ or $B_2$. Then the
subgroup $[\aut (L),G(L)]$ of the symmetric group on the
elements of $L$ is isomorphic to the right nucleus of $L$. 
\end{prop}
\bp Take an element $\alpha \in \aut (L)$ and a generating
element $\lambda_x^{-1}$ of $G(L)$. We have
\[ [\alpha, \lambda_x^{-1}] = \alpha^{-1} \lambda_x \alpha
\lambda_x^{-1} = \lambda_u \> \lambda_u^{-1} \lambda_{x^\alpha}
\lambda_x^{-1} = \lambda_u \tau(x,\alpha),\]
where $u \in L$ such that $x^\alpha u^{-1} = x$ and
$\tau(x,\alpha) = \lambda_u^{-1} \lambda_{x^\alpha}
\lambda_x^{-1}$. 

Using again $[L : N_\rho] = 2$, we obtain that $u \in N_\rho$.
Furthermore, we know that any element of $N_\rho$ can occur
as $u$, because the group of automorphisms acts transitively on
the set $L\setminus N_\rho$. By the remark after the \ref{stab},
$\lambda_u$  commutes with the involution $\tau(x,\alpha)$.
\ref{stab} says that $\tau$ depends only on $u$, since
$\tau(x,\alpha) = id$ if and only if among $u, x, x^\alpha$
there is a left nucleus element iff $u$ is a left nucleus
element. But the left 
nucleus $N_\lambda$ has index 2 in the right nucleus $N_\rho$,
and the nontrivial $\tau(u)$ has order 2, thus the map
$[\alpha, \lambda_x^{-1}] = \lambda_u^{-1} \tau(u) \mapsto u$
define an isomorphism between the set of generators of
$[\aut (L), G(L)]$ and the group $N_\rho$. 
\ep
%
%
%
%
\section{Direction preserving collineations}

In this and the next sections, we are going to determine the
group of direction preserving collineations in the smallest Bol
3-nets.
For that, let us fix some notations. In the following, $\Gamma(\cal
N)$ will always denote the direction preserving collineation
group of the Bol 3-net $\cal N$. If it does not cause
misunderstanding, we only write $\Gamma$. If the net is 
coordinatized by the loop $(L,\cdot)$, then $\Gamma_0({\cal N})$
is going to be the subgroup stabilizing the origin. By $l_h$ and
$l_v$ we denote the horizontal and vertical lines through the
origin, respectively, and $\Gamma_H({\cal N})$ and
$\Gamma_V({\cal N})$ are the stabilizers of $l_h$ and $l_v$ in
$\Gamma({\cal N})$, respectively. We use the notation $P$ for
the orbit of the origin in the action of $\Gamma({\cal N})$ on
the point set of $\cal N$. Furthermore, we introduce the notation
$H = N({\cal N}) \Gamma_0({\cal N}) \leq \Gamma(\cal{N})$.

Now, we state our Main Theorem, whose proof needs this last
chapter. 
\begin{mth}
Let the $3$-net $\cal N$ coordinatized by $B_1$ or $B_2$. Then,
using the above introduced notations, we can state the following.
\setcounter{szam}{0}
\begin{list}{(\roman{szam})}{\usecounter{szam}}
\item In both cases $B_1$ and $B_2$, we have $H^\prime = \Phi
(H)$ and these groups are abelian.
\item For $B_1$,
\[\Gamma = H = \Gamma_0 \times_\varphi N \cong \aut (B_1)
\times_\varphi G(B_1).\]
\item For $B_2$,
\[\Gamma = H \times C_2,\]
moreover there exists an Abelian normal subgroup $\Lambda$ of
$G(B_2)$, which operates sharply transitively on the set $B_2$,
and holds
\[H \cong Aut(B_2) \times_\varphi \Lambda.\]
\end{list} \end{mth}
First, consider the kernel of the map $\Phi: N({\cal N})
\rightarrow G(L)$, where the loop $(L,\cdot)$ is equal to $B_1$
or $B_2$ and the $3$-net $\cal N$ is coordinatized by
$(L,\cdot)$. With some calculation, one can obtain that in the
loop $B_1$, the map $J:x\mapsto x^{-1}$ is a right
pseudo-auto\-mor\-phism but not an automorphism. On the other
hand in the loop $B_2$, $J$ is an automorphism, more precisely
$J = \psi_1^2 \circ \psi_2$. By \ref{kernel}, this means,
that in the smallest Bol $3$-nets, the kernel of $\Phi$ is
trivial. The elements of $H$ are collineations, this means they
are pairs of two permutations, and the projection of the second
permutation  of these col\-linea\-tions gives us an isomorphism
from $H$ onto the subgroup 
\[\bar{H} = G(L) \aut (L)\]
of the symmetric group on the
carrier set $L$. From now on, we will consider the group
$\bar{H}$. 

It is easy to see that the commutator subgroup $\bar{H}^\prime$
is generated by ${\aut (L)}^\prime$ and $[G(L), \aut (L)]$,
$G(L)^\prime$. 

In the paper \cite{NSt}, the generated groups $G(B_1)$ and
$G(B_2)$ are determined. Their result was
that they are non-abelian groups of order 16, such that their
commutator subgroups have order 2 and contain the square
elements of the group. Actually, also the dihedral group $D_8$
of order 8 has this last property, its commutator subgroup
contains the squares of the group. 
\begin{theo}
Let the loop $L$ be equal to $B_1$ or $B_2$. Then the commutator
subgroup $\bar{H}^\prime$ is isomorphic to $C_2 \times C_4$ or
$C_2 \times C_2 \times C_2$, respectively. 
\end{theo}
\bp Firstly, we show that $G(L)^\prime \subseteq [\aut (L),
G(L)]$. By \cite{NSt} we know that the group $G(L)^\prime$ is a 
group of order 2, and its non-trivial element is the left
translation by the non-trivial element $n_0$ of the left
nucleus. On the other hand, using the notation of \ref{nuc}, we
can choose $u = n_0$, and so we have $u \in N_\lambda$, 
so $\tau(u) = id$, thus $\lambda_{n_0} \in [\aut (L), G(L)]$. 

In the second step, we show that $[\aut (L), G(L)] \cap
{\aut (L)}^\prime = \{ id \}$. As we know, the elements of
$[\aut (L), G(L)]$ are in the form $\lambda_u^{-1} \tau(u)$,
with $u\in N_\rho$. This permutation fixes the unit element $1$
of the loop $L$ only if $u=1$, but any element of ${\aut
(L)}^\prime$ must fix $1$.

Finally, we show that elements of ${\aut (L)}^\prime$ commute
with elements of the group $[\aut (L), G(L)]$. In both cases,
the group 
$\aut (L)$ of automorphisms is isomorphic to the dihedral group
$D_8$ of order 8. The commutator of this group contains that
unique involution which commutes with any other element.
On the other hand, we saw already, that $\lambda_a \lambda_b
\lambda_c$, with $c\cdot ba =1$ is an involutory automorphism.
It can take two different values, the identity or the
non-trivial $\tau (u)$. One can easily check that, for both
cases, it commutes with any other automorphism of the loop $L$.
Thus, the non-trivial $\tau (u)$ is the only non-trivial element
of the commutator group ${\aut (L)}^\prime$. But we had already
noticed that this element fixes the elements of the right
nucleus, which means that it commutes with an element of the
form $\lambda_u^{-1} \tau (u)$.

So, the theorem is proven, hence we have shown that
\[\bar{H}^\prime \cong C_2 \times N_\rho,\]
and $N_\rho \cong C_4$ in the loop $B_1$ and $N_\rho \cong C_2
\times C_2$ in the loop $B_2$. 
\ep
The nice structure of these smallest Bol loops gives us another
interesting result. 
\begin{theo}
If the loop $L$ is equal to $B_1$ or $B_2$ and $\bar{H}$ is the
group $\bar{H} = \aut (L) G(L)$ as before, then the Frattini
subgroup of $\bar{H}$ is the commutator subgroup
$\bar{H}^\prime$.
\end{theo}
\bp It is known that for a $p$-group, the Frattini subgroup
$\Phi (\bar{H})$ is the smallest normal subgroup $K$ such that
$\bar{H} / K$ is elementary abelian, and so it must contain the
commutator group. Thus, we only have to show that the square of
any element of $\bar{H}$ is contained in the commutator subgroup
$\bar{H}^\prime$. 

Any element of $G(L)$ is in the form $\tau
\lambda_x$, where $\tau \in G(L)$ stabilizes the unit element
$1$ and so it is an automorphism of $L$. This means that an
element of $\aut (L) G(L)$ is of the form $\beta \tau \lambda_x =
\alpha \lambda_x$, where $\beta \in \aut (L)$ and $\alpha =
\beta \tau \in \aut (L)$.
Then,
\[\alpha \lambda_x \alpha \lambda_x = \alpha^2
\lambda_{x^\alpha} \lambda_x = 
\alpha^2 \lambda_u^{-1} \lambda_u \lambda_{x^\alpha}
\lambda_x^{-1} \lambda_{x^2} =\]
\[= \alpha^2 \lambda_u^{-1} \tau(u) \lambda_{x^2} \in (\aut
(L))^\prime [G(L), \aut (L)] G(L)^\prime = (\aut (L)
G(L))^\prime.\] 
Thus, we obtain that the Frattini subgroup $\Phi (\bar{H})$ is
contained in the commutator subgroup $\bar{H}^\prime$, and so
they must be equal. \ep
%
%
%
%
In the following two subsections, we will consider the 3-nets
determined by the loops $B_1$ and $B_2$. We will collect
information about their collineation group, and for that we use
two methods which are different but rather similar. Anyway, both
methods can be a useful tool in the examination of collineation
groups of more general Bol 3-nets.
\subsection{Collineations of the 3-net $B_1$}
\begin{prop} \label{regu}
Let $\cal N$ be a Bol 3-net coordinatized by the Bol loop
$(L,\cdot)$. Assume that $J$ is a right pseudo-auto\-mor\-phism
of $L$ but not an automorphism and that the right nucleus has
index $2$ in the loop $L$. Then $N=N(\cal N)$ acts semi-regularly
on $P$, the orbit of the origin in $\cal N$.
\end{prop}
\bp Suppose that $(\alpha,\beta)\in N$ and $1^\alpha = 1^\beta = 1$.
The permutation $\beta \in G(L)$ fixes the unit element $1$,
thus it is generated by elements of the form $\lambda_a
\lambda_b \lambda_c$, with $c\cdot ba = 1$. Then $\beta$ fixes
the elements of the right nucleus $N_\rho$, since they are fixed
by the elements of the type $\lambda_a \lambda_b \lambda_c$, as we
have seen it in the remark after \ref{stab}. Another
remark, after \ref{kernel},  says that if $\beta$ fixes
the unit element $1$ then $\alpha$ fixes those elements which
can occur as the right companion of $J$. By the conditions of this
proposition, this means that $\alpha$ fixes any element outside
of $N_\rho$. So $(\alpha, \beta)$ fixes the origin, hence
$\alpha = \beta$, thus $\alpha = \beta = id$.
\ep
The loop $B_1$ is left conjugacy closed (cf. \cite{Burn}), so
the orbit $P$ is a 
union of vertical lines. Moreover, we have a bijection from $l_h
\cap P$ onto $N_\lambda$, and so $|P|=16$. Comparing this with
\ref{regu}, the result is that $N$ acts regularly on
$P$. This means $\Gamma = \Gamma_0 N$ and $\Gamma_0 \cap N =
\{1\}$. Applying the isomorphy theorem of groups,
\[\Gamma/N = \Gamma_0 N/N \cong \Gamma_0/{\Gamma_0 \cap N} =
\Gamma_0.\]
So we got that $\Gamma$ splits over $N$, with complement
$\Gamma_0$. From the group theory it is known that in this case we
have the semidirect product
$\Gamma = \Gamma_0 \times_\varphi N$, where the action 
$\varphi$ of $\Gamma_0$ on $N$ is defined by restriction of the
action of $\Gamma$ on $N$ by conjugation. 
\begin{theo}
Let $\cal N$ be the Bol 3-net coordinatized by the Bol loop
$B_1$ and let $\Gamma$ be the group of the direction preserving
collineations of $\cal N$. Moreover, denote by $\varphi$ the
action of $\aut (B_1)$ on $G(B_1)$, defined by $\alpha :
\lambda_x \mapsto \lambda_{x^\alpha}$. Then
\[\Gamma \cong \aut (B_1) \times_\varphi G(B_1).\]
\end{theo}
\bp We know already that $\Gamma$ is a semi-direct product of
this form, $\alpha^{-1} \lambda_x \alpha = \lambda_{x^\alpha}$.
If we conjugate inside of the collineation 
group, the action will be the same, since the isomorphisms from
the subgroups of $\Gamma$ onto $G(B_1)$ and $\aut (B_1)$ are
projections in fact. \ep
\begin{prop}
The commutator subgroup $\Gamma^\prime$ is equal to the Frattini
subgroup $\Phi (\Gamma)$ and they are isomorphic to the
abelian group $C_2 \times C_4$.
\end{prop}
\bp If we denote $H = \Gamma_0 N$, then by the last theorem we
have $\Gamma = H \cong \bar{H}$. We have already proven that
$\Gamma^\prime \cong \bar{H}^\prime = \Phi (\bar{H}) \cong \Phi
(\Gamma)$ is isomorphic to $C_2 \times C_4$. \ep
\subsection{Collineations of the 3-net $B_2$}
Unfortunately, in the case of $B_2$ we can not apply the above
method, because the map $J$ is an automorphism, and so $N(\cal
N)$ does not act semi-regularly on $P$, the orbit of the origin.
\begin{prop} \label{n0}
Let $\cal N$ be a Bol 3-net coordinatized by the Bol loop
$(L,\cdot)$ and suppose that $J$ is an auto\-mor\-phism of $L$.
Then 
\[N_0 \cong G_1 < \aut (L),\]
where $N_0$ is the stabilizer of the origin in the normal
subgroup $N(\cal N)$ and $G_1 < G(L)$ contains the elements
fixing the unit element $1$ of the loop.
\end{prop}
\bp We saw already, that if we coordinatize the 3-net $\cal N$
by $(L,\cdot)$, then the elements of $N(\cal N)$ are in the
form 
\[(\rho_u \lambda_{m_1} \cdots \lambda_{m_k}, \lambda_{m_1}^{-1}
\cdots \lambda_{m_k}^{-1}) = (\alpha, \beta),\]
with $u = m_1(\cdots(m_{k-2}\cdot m_{k-1} m_k)\cdots ) =
1^{\beta^{-1}}$. If the map $J$ is an
automorphism, then 
\[1^{\lambda_{m_1} \cdots \lambda_{m_k}} =
(1^{\lambda_{m_1}^{-1} \cdots \lambda_{m_k}^{-1}})^{-1},\] 
That means that in our case,
\[\beta \in G_1 \Leftrightarrow u=1 \mbox{ and } 1^\alpha =
(1^\beta)^{-1} \Leftrightarrow (\alpha,\beta) \mbox{ fixes }
1,\] 
and so we have a bijection from $G_1$ onto $N_0$, which is the
restriction of $\Phi^{-1}$, thus, it is an isomorphism. \ep

By [4], we have $|G_1|\cdot |B_2| = |G(B_2)| = 16$, hence $G_1 =
\{id, \tau \}$, where $\tau = \lambda_a \lambda _b \lambda_c$
with $c \cdot ba = 1$. By \ref{n0}, 
\[N_0 = \{ (id, id), (\tau,\tau)\} = N \cap \Gamma_0 \lhd
\Gamma_0.\] 

Now, let us consider the collineation group of the 3-net as a
permutation group on the horizontal lines. From the paper \cite{Funk},
we know that the subgroup $N=N(\cal N)$ acts transitively on this
set of lines. By \ref{dircoll}, an element of the stabiliser of
such a horizontal line, say $l_h$, determines a left
pseudo-auto\-mor\-phism of the loop $B_2$, with a companion
element. However, any left pseudo-auto\-mor\-phism of $B_2$ is
an automorphism and the left companion elements form the group
$N_\lambda$. Thus, the stabiliser of $l_h$ is isomorphic to
the semidirect product $\aut (B_2) \times_\varphi N_\lambda(B_2)$,
since an automorphism of the loop induces an automorphism on the
left nucleus (see \ref{pseu}). Moreover, the
left nucleus of $B_2$ has 2 elements, so its only automorphism
is the identity, hence we have
\begin{prop}
Let $\cal N$ denote the $3$-net coordinatized by the Bol loop
$B_2$. Then for the stabiliser $\Gamma_H$ of the vertical line
$l_h$ through the origin in the direction preserving
collineation group of the net, we have the direct product
\[\Gamma_H \cong D_8 \times C_2.\]
\nop \end{prop}

\begin{prop}
For the 3-net $\cal N$, coordinatized by the Bol loop $B_2$, we
have 
\[\Gamma / N \cong (\aut (B_2) / G_1) \times N_\lambda \cong C_2
\times C_2 \times C_2.\]
\end{prop}
\bp By the isomorphy theorem of groups,
\[\Gamma/N = \Gamma_H N/N \cong \Gamma_H/(\Gamma_H \cap N).\]
Obviously, $N_0 \subseteq \Gamma_H \cap N$, since if a
collineation fixes the origin, it fixes the line through it, as
well. On the other hand, if $(\alpha, \beta) \in \Gamma_H \cap
N$, then $1^\beta = 1$, and by \ref{n0}, $1^\alpha
= 1$, thus $N_0 = \Gamma_H \cap N$. But $N_0$ is normal in
$\Gamma_0$ which is normal in $\Gamma_H= \Gamma_0 \times C_2$,
so 
\[\Gamma_H/(\Gamma_H \cap N) = (\Gamma_0 \times C_2)/N_0 =
(\Gamma_0/N_0) \times C_2 \cong (\aut (B_2) / G_1) \times
N_\lambda.\] 
We still have to show that $\aut (B_2) / G_1 \cong C_2 \times
C_2$. The group $\aut (B_2)$ is isomorphic to $D_8$, which
contains 5 involutions, but only one of them commutes with all
the others. In the case of $\aut (B_2)$ this is $\tau = \lambda_a
\lambda_b \lambda_c \in G(L)$, with $c \cdot ba = 1$. Now, if we
factorize by this special involution, we obtain exactly the
factor group $C_2 \times C_2$. \ep
\begin{coro}
Denote the subgroup $\Gamma_0 N$ of the collineation group
$\Gamma$ by $H$ and by $\sigma$ the collineation
$(\lambda_0,id)$ where $\lambda_0$ is the left translation with
the nontrivial element of the left nucleus of $B_2$. Then
\setcounter{szam}{0}
\begin{list}{(\roman{szam})}{\usecounter{szam}}
\item $\Gamma = H \times \{ (id,id), \sigma \}$;
\item the commutator subgroup $\Gamma^\prime$ is isomorphic to
the abelian group $C_2 \times C_2 \times C_2$;
\item the Frattini subgroup of $\Gamma$ is equal to the
commutator subgroup, thus, $\Gamma$ is generated by $4$ elements.
\end{list} \end{coro}
\bp Obviously, $\sigma$ is an involutory collineation of the
3-net. From the definition of the left nucleus elements follows
that $\lambda_0$ commutes with all right translations. By [3],
$\lambda_0$ is contained in the center of the 
group $G(B_2)$, and so $\lambda_0$ commutes with any
element of $N$. The left nucleus has order 2, an automorphism of
the loop induces an automorphism of the nuclei, hence the
elements of the left nucleus are fixed by every automorphism.
With other words, $\sigma = (\lambda_0,id)$ commutes with each
element of $\Gamma_0$. 

Moreover, if an element of $N$ fixes the
horizontal line $l_h$ through the origin, then it fixes the
origin, too. So does any element $\gamma_0 \nu \in H$, where
$\gamma_0 \in \Gamma_0$ and $\nu \in N$, hence $\sigma$ is not
contained in $H$. And the subgroup $H$ has index 2, since for
the factor groups $|H/N| = 4$ and $|\Gamma /N| = 8$, thus (i) is
proven.

From this it follows that $\Gamma^\prime = H^\prime \cong C_2
\times C_2 \times C_2$. 

Clearly, the Frattini subgroup $\Phi (\Gamma)$ is generated by
the squares in $\Gamma$, this means by (i) that $\Phi (\Gamma) =
\Phi (H) = H^\prime$, see Theorem 5.5. Finally, we have 
\[\Gamma / \Gamma^\prime \cong C_2 \times C_2 \times C_2 \times
C_2,\] 
the factor group has order $2^4$, but than the minimal number of
element which generate $\Gamma$ is 4 (cf \cite{Kurz}, p. 58.).
\ep
\begin{theo}
There exists an Abelian normal subgroup $\Lambda$ of
$G(B_2)$, which operates sharply transitively on the set $B_1$,
and holds $H \cong Aut(B_2) \times_\varphi \Lambda$.
\end{theo}
\bp By \cite{NSt}, we know that there is a unique Abelian normal
subgroup $\Lambda$ of $G(B_2)$, which operates sharply
transitively on $B_2$. To this $\Lambda$, a Abelian normal
subgroup $N_\Lambda$ of
$N$ corresponds, which acts sharply transitively on the set of
vertical lines. Its uniqueness means that it is a normal
subgroup in $H$ as well. Then we have that $H / N_\Lambda$ is
isomorphic to the stabilizer of the vertical line through the
origin in $H$, which is equal to the stabilizer of the origin in
$\Gamma$, thus we have the desired semidirect product $H \cong
Aut(B_2) \times_\varphi \Lambda$. 
\ep 
And by this, the proof of the Main Theorem is complete. Finally,
let us turn back to the full collineation group. We saw already,
that any collineation of a non-Moufang Bol $3$-net leaves the
set of vertical lines fixed. This means that a collineation
which does not preserve the directions must interchange the
horizontal and the transversal directions, the product of two
such ones is a direction preserving collineation. Thus, the
$\Gamma(\cal N)$ is a normal subgroup of index 2 in the full
group of collineations of the $3$-net $\cal N$.
%
%
%
%

G\'abor Nagy

Bolyai Institute

J\'ozsef Attila University

Aradi v\'ertan\'uk tere 1.

6720 Szeged (Hungary)
\end{document}